\documentclass [a4paper,twoside,12pt]{article}

\usepackage[utf8]{inputenc}

\usepackage[T1]{fontenc}
\usepackage{amsmath,amsfonts,amssymb}
\usepackage{vmargin,graphicx,theorem}
\usepackage[french]{babel}
\usepackage{enumerate}
\usepackage{color}
\usepackage{pst-fill,pst-grad,pst-plot,pst-eucl,pstricks-add,pst-node}

\setpapersize[portrait]{A4}
\setmarginsrb{1.5cm}{1cm}{1.5cm}{2.5cm}{1.5cm}{0cm}{0.5cm}{2cm}
\newcommand{\disp}{\displaystyle}
\newcommand{\dR}{\ensuremath{\mathbf{R}}}

\newcommand{\proofend}{~$\rhd$}
\newcommand{\proofbegin}{~$\lhd$}

\newcommand{\PAR}[1]{\ensuremath{{\left(#1\right)}}} % (1)
\newcommand{\SBRA}[1]{\ensuremath{{\left[#1\right]}}} % [1]
\renewcommand{\phi}{\varphi}
\renewcommand{\geq}{\geqslant}

\def\disp{\displaystyle}

\newcommand{\R}{\dR}

\newcommand{\beq}{\begin{equation}}\newcommand{\eeq}{\end{equation}}

\begin{document}

\title{L'entropie, de Clausius aux inégalités fonctionnelles}
\author{ Ivan Gentil}

\date{\today}

\maketitle

Nous nous intéressons dans ce document  à  l'entropie~\footnote[1]{Ce document  a été préparé à l'occasion d'un colloquium donné le 12 avril 2019 à l'Université Paris Descartes que je  remercie pour l'invitation.}. L'entropie est multiple, l'idée est de la décrire  dans le prolongement de sa définition proposée par le physicien Clausius. En effet, Clausius expose en 1865 le second principe de la thermodynamique et  propose aussi le concept d'entropie. Au lieu de définir simplement une fonctionnelle, point central pour l'élaboration du second principe, il va en fait définir un concept suffisamment  général  pour qu'il soit utilisé dans de  nombreux domaines des mathématiques.

Dans ces quelques pages, je souhaite faire apparaître le rôle joué par l'entropie dans le domaine de l'étude des flots de gradient  et des inégalités fonctionnelles. 
Partant de la définition de Clausius en 1865, je vais tenter d'expliquer  comment des inégalités fonctionnelles fondamentales comme l'inégalité de Sobolev, point clé en analyse, sont des inégalités de structures naturelles reliées à une certaine entropie.  Ce cheminement me permet de donner un aperçu de l'utilisation des flots de gradient en dimension finie ou infinie,   de la théorie de Bakry-\'Emery et plus récemment du calcul d'Otto.

\section{L'entropie proposée par Clausius}

% Dans Zentrablatt le mot {\it Entropy} apparait  41000 fois et 480000 pour le mot {\it function}.  Notons que Zentralblatt a 3 millions d'entrées depuis 1868 soit en gros à l'apparition du mot entropie dans l'article de  Clausius. Si un mot est référencé c'est qu'il est dans le titre, dans le résumé ou bien dans les références.  Ainsi c'est un peu plus de 1\% des entrées.  Pour le mot  {\it Probability} c'est 200 000 entrées. 
On retrouve 6 millions de fois le terme Entropy sur le moteur de recherche Ecosia et près de 40 millions sur Google. Il apparaît 41000 fois dans Zentrablatt, ce qui représente un peu plus de 1\% des entrées.  Ce terme est tellement utilisé que l'on ne peut pas être exhaustif dans une courte introduction, il en existe d'ailleurs bien d'autres.

Une fois n'est pas coutume, on sait quand ce mot a été inventé, c'est en 1865 dans un article de  Rudolf Julius Emmanuel Clausius, physicien allemand du milieu du XIXe siècle, cf.~\cite{clausius}.

Avant de parler du contenu scientifique de l'article, arrêtons-nous sur sa forme. Cet  article est publié dans la revue de mathématiques  {\it Journal de mathématiques pures et appliquées} en 1865.  Rappelons que c'est un des plus vieux journaux de mathématiques existants, fondé par 
Joseph Liouville en 1836, juste après le  journal de Crelle (Journal für die reine und angewandte Mathematik) qui a été fondé, lui, en 1826. On peut se demander 
pourquoi Clausius choisit un journal de mathématiques  et pourquoi celui-ci en langue française. Le journal de Crelle existe et est allemand~!  Est-ce le prestige de ce nouveau journal 
français ? C'est 
d'autant plus étonnant car Clausius a publié la très grande majorité de ses articles en allemand, dans {\it Annalen der Physik }, dont la publication commence dès
1790 et aussi dans le journal de Crelle.  Ce qui est certain c'est qu'à cette époque, le choix du journal est fondamental pour une diffusion rapide des idées contrairement 
 à maintenant, où le choix d'une 
 revue est plus lié à son prestige et à la carrière du chercheur.  Clausius n'écrivait probablement pas assez bien le français puisque son article a été traduit par un 
traducteur qui semble professionnel et non mathématicien. 
%\medskip 

Cet article  est clairement important car Clausius y  formule le second principe de la thermodynamique en définissant l'entropie.  On peut  en dégager   deux points importants. 
\begin{enumerate}[\bf 1.]
\item  Il participe à l'élaboration du second principe de la thermodynamique~: l'irréversibilité des phénomènes physiques lors d'échanges thermiques. ll s'appuie  sur des travaux précédents du physicien Sadi Carnot, mort à 36 ans en 1832. Notons que Carnot n'a publié qu'un seul ouvrage en 1824 intitulé {\it Réflexions sur la puissance motrice du feu}, livre  fondateur de la thermodynamique, même si ce mot n'a été inventé que plus tard. Carnot souhaite, dans son livre,  améliorer ce qu'il appelle la machine à feu, dont l'exemple fondamental est  la machine à vapeur développée à la fin du $\rm XVIII^e$.
\item Pour démontrer ce second principe Clausius introduit le concept d'entropie~: exhiber une fonctionnelle qui montrera l'irréversibilité d'un phénomène physique. Il comprend qu'il y a deux notions importantes, l'énergie et ce nouveau concept. Il a donc cherché un mot proche du mot {\it énergie}. Son idée est d'utiliser la racine du mot grecque $\eta\tau\rho o\pi\eta$ qui signifie transformation, ou plutôt {\it en transformation}. L'utilisation de cette racine grecque lui donne une consonance universelle et le mot \og entropie\fg  s'écrit presque de la même façon dans toute les langues utilisant l'alphabet latin. Il pense, et l'histoire lui a donné raison, que ce mot représente un concept général. 

Ce mot va effectivement rester  et être largement utilisé dans de nombreux domaines. 
%C'est une façon intelligente pour que ce concept se propage auprès de la communauté  scientifique. 
Clausius note l'entropie  $S$ et l'assimile à {\it du désordre qui ne peut que croître}.  Citons l'exemple récent et emblématique de Perelman qui définit aussi sa propre entropie, fonctionnelle centrale dans sa preuve de la conjecture de Poincaré.  Dans sa célèbre prépublication~\cite{perelman} de 2002, Perelman note aussi l'entropie $S$, comme Clausius. 
%Bien entendu Perelman n'est pas le seul à sa propre entropie. 

%, il n'y a pas une enparticulière mais plutôt un concept. 

\end{enumerate}
 
\medskip 

Clausius  termine son papier par des considérations générales, bien connues maintenant,  qui sont les deux principes de la thermodynamique. Nous les reproduisons tels qu'ils sont écrits dans son article page 400~: 
\begin{itemize}
\item \it L'énergie de l'univers est constante
\item L'entropie de l'univers tend vers un maximum
\end{itemize}
A l'école, j'ai plutôt appris  que   {\it l'entropie d'un système isolé est croissante}. 

\bigskip

Illustrons simplement le concept proposé par Clausius avec l'équation de la chaleur dans $\dR^n$. Partant d'une mesure de probabilités  $\mu_0$ dans $\R^n$, c'est-à-dire une mesure positive de masse 1, l'équation de la chaleur est l'équation aux dérivées partielles suivante, 
$$
\partial_t \mu_t=\Delta \mu_t,\, t>0.
$$
Dans le cadre de cet article, on peut allègrement remplacer les mesures de probabilité par des  fonctions positives ayant une intégrale par rapport à la mesure de Lebesgue égale à 1.  Rappelons que partant d'une mesure de probabilités $\mu_0$ sur $\R^n$, l'équation de la chaleur admet une unique solution régulière sur $\R_{+}^*\times\R^n$, qui est un flot dans l'ensemble des mesures de probabilité dans $\R^n$, noté $\mathcal P(\R^n)$, 
$$
[0,\infty)\ni t\mapsto \mu_t\in\mathcal P(\R^n).
$$ 
Dès 1872, soit seulement quelques années après l'article de Clausius, Boltzmann a compris les deux points fondamentaux de Clausius~: {\it l'irréversibilité  du phénomène physique} et le {\it concept de l'entropie}. Ainsi il définit sa propre entropie, appelée maintenant entropie de Boltzmann. Pour toute mesure de probabilité $\mu\in\mathcal P (\R^n)$ admettant une densité (que l'on note abusivement aussi $\mu$) par rapport à la mesure de Lebesgue $\lambda$,   
\begin{equation}
\label{eq-9}
{\rm Ent}(\mu)=\int \mu\log\mu=\int \log\frac{d\mu}{d\lambda}\, d\mu\in[0,\infty].
\end{equation}
Cette entropie est proposée par Boltzmann en 1872 dans son célèbre théorème 
$H$, (nous proposons ici seulement une version simplifiée).  L'entropie de Boltzmann diffère de l'entropie thermodynamique de Clausius mais c'est la 
même que celle utilisée par Shannon au $\rm XX^e$ siècle, avec un signe opposé.  Citons aussi le cas de Nash qui utilise aussi l'entropie de Boltzmann en 1958 pour démontrer, dans un article fondateur et de façon surprenante, la régularité d'équations elliptiques et paraboliques. 
%Cet article est précurseur car personne auparavant n'avait pu imaginer cette utilisation.  

Ainsi, pour simplifier, partant d'une donnée initiale $\mu_0$ à densité régulière (par exemple dans l'espace de Schwartz), si on dérive l'entropie de Boltzmann le long du flot de la chaleur,  on obtient pour $t\geq0$,
\begin{equation}
\label{eq-13}
\frac{d}{dt}{\rm Ent}(\mu_t)=\int (1+\log \mu_t)\partial_t \mu_t=\int \Delta \mu_t\log \mu_t=-\int \nabla \mu_t\cdot\nabla\log\mu_t=-\int \frac{|\nabla\mu_t|^2}{\mu_t},
\end{equation}
où le point clé est l'utilisation d'une intégration par parties dans la troisième égalité.  
L'entropie du système est toujours décroissante. On ne peut pas revenir en arrière dans une évolution naturelle. Si la condition initiale est par exemple  une masse de Dirac en 0 ($\mu_0=\delta_0$ et dans ce cas le calcul donné en~\eqref{eq-13} n'est valable que pour $t>0$), la chaleur va se diffuser dans l'espace tout entier. Il est impossible de revenir à une masse de Dirac, autrement dit la diffusion de la chaleur a un sens, celui donné par {\it la flèche du temps}. Nous concluons que {\bf l'évolution de cette transformation physique est irréversible.} Il est important de remarquer que pour démontrer l'irréversibilité de l'équation de la chaleur, la fonctionnelle d'entropie est seulement un outil. L'entropie de Boltzmann est une fonctionnelle remarquable, comme nous le verrons en section~\ref{sec-3}, mais il y en d'autres qui montrent d'irréversibilité de l'équation de la chaleur, comme par exemple la norme $L^p$ ($p>1$) de la densité, $\mu\mapsto \int \mu^p$. 

%\medskip

Bien entendu, toutes les transformations physiques ne sont pas forcément irréversibles. Par exemple, la chute libre d'un corps  satisfait à l'équation de Newton qui est,  elle, une équation  réversible. Si on renvoie  le corps avec une vitesse inversée, il retournera à sa place initiale~! 
%L'équation de Newton est une équation  d'ordre~2 alors que l'équation de la chaleur peut-être vue comme une équation d'ordre~1, ce que nous verrons en section~\ref{sec-3}.

%\medskip 

%Il est amusant et anecdotique de noter que les terminologies sont assez paradoxales. Si comme nous venons de le voir,  l'équation de la chaleur est une équation irréversible, elle repose  sur le mouvement brownien qui est lui un processus stochastique réversible, si on renverse le temps (partant ici de la mesure de Lebesgue) la loi de l'ensemble du  processus ne change pas. 

\medskip

En conclusion, le concept d'entropie est un outil qui s'adapte à ce que l'on veut étudier. Nous allons illustrer, dans les sections suivantes,   son utilisation remarquable dans l'étude des flots de gradient, en dimension finie qui est un cadre simple pour énoncer les propriétés et en dimension infinie, où l'on retrouve des résultats importants d'analyse fonctionnelle.

\section{Flots de gradient en dimension finie et entropie}
\label{sec-2}

Nous explorerons maintenant le cas simple d'un flot de gradient, permettant d'illustrer la théorie de Bakry-\'Emery proposée en  1985 dans~\cite{bakry-emery1985}.
\medskip

Considérons une fonction de classe $\mathcal C^2$,  
$$
E:\dR^n\mapsto\dR,
$$ 
vérifiant pour un certain $\rho>0$, la condition   
\begin{equation}
\label{eq-2}
{\rm Hess\, }E\geq \rho {\rm Id},
\end{equation}
où l'inégalité est vue dans le sens des matrices symétriques. 
Cette fonction est alors coercive, elle vérifie $\lim_{|x|\rightarrow\infty}E(x)=\infty$, et admet donc un minimum global atteint en un unique point $\beta\in\R^n$. On garde en mémoire l'exemple typique,  $E(x)={|x|^2}/{2}$ où $\rho=1$ et $\beta=0$.

\medskip

$\bullet$ {\bf Flot de gradient} 

On note  $(S_t(x))_{t\geq 0}$, $x\in\dR^n$,  la solution de l'équation différentielle ordinaire 
$$
\left\{
\begin{array}{l}
\dot{X}_t=-\nabla E(X_t);\\
X_0=x .
\end{array}
\right.
$$
On utilise la notation provenant de la mécanique $\dot{X}_t=\frac{d}{dt}X_t$, désignant le vecteur  vitesse de la trajectoire $(X_t)_{t\geq0}$. 

On dit  que $(S_t)_{t\geq 0}$ est {\it le flot de gradient de $E$ par rapport à la métrique euclidienne, celle qui définit le gradient usuel}. 
Dans l'exemple  classique, $E(x)=|x|^2/2$, la solution partant de $x$ est simplement $S_t(x)=e^{-t}x$, solution d'une équation différentielle ordinaire, linéaire et de degré~1.

\medskip

$\bullet$ {\bf Identité de  de Bruijn} 

Nous appelons $E$ la fonctionnelle d'entropie du système, on a  
\begin{equation}
\label{eq-5}
\frac{d}{dt}E(S_t(x))=\nabla E(S_t(x))\cdot\frac{d}{dt}S_t(x)=-|\nabla E(S_t(x))|^2\leq0.
\end{equation}
L'opposé du terme de droite est appelé {\it la production d'entropie} (appelé parfois énergie, peut-être  à tort) et cette identité porte parfois  le nom de de Bruijn.

Pour les mêmes raisons que précédemment, ce phénomène physique est irréversible. Si on imagine $S_t(x)$ comme étant la position d'une bille à l'instant $t$, celle-ci va descendre vers $\beta$, le point où $E$ réalise son unique minimum global. C'est exactement la même chose que la chaleur, la bille ne  remontera pas la pente.

\medskip

$\bullet$ {\bf Méthode de Bakry-\'Emery} 

Clausius invente le concept d'entropie, Boltzmann propose de dériver l'entropie le long du flot.  
L'idée fondamentale de Bakry et \'Emery est de dériver une seconde  fois l'entropie le long du flot, 
\begin{multline*}
\frac{d^2}{dt^2} E(S_t(x))=-\frac{d}{dt}|\nabla E(S_t(x))|^2=-2\nabla E(S_t(x))\cdot\frac{d}{dt}\nabla E(S_t(x))=\\
2\nabla E(S_t(x))\cdot{\rm Hess\, }E(S_t(x))\nabla E(S_t(x)).
\end{multline*}
En utilisant l'inégalité ${\rm Hess\, }E\geq \rho {\rm Id}$, on obtient 
$$
\frac{d}{dt}|\nabla E(S_t(x))|^2\leq -2\rho|\nabla E(S_t(x))|^2,
$$
soit donc après intégration sur l'intervalle $[0,t]$, 
\begin{equation}
\label{eq-10}
|\nabla E(S_t(x))|^2\leq e^{-2\rho t}|\nabla E(x)|^2.
\end{equation}
Sous la condition de convexité~\eqref{eq-2}, la production d'entropie décroit vers 0 avec une vitesse  exponentielle et explicite. 

\medskip

$\bullet$ {\bf Inégalité  entropie-production d'entropie} 

On peut montrer que $\lim_{t\rightarrow\infty}S_t(x)=\beta$. En effet, $\nabla E(S_t(x))$ tend vers 0  grâce à l'inégalité~\eqref{eq-10} et $\beta$ est le seul point qui annule le gradient de $E$.   

Ainsi on a 
\begin{multline*}
E(x)-E(\beta)=E(S_0(x))-E(\lim_{t\rightarrow\infty}S_t(x))=-\int_0^\infty \frac{d}{dt}E(S_t(x))dt=\\
\int_0^\infty|\nabla E(S_t(x))|^2dt\leq \int_0^\infty e^{-2\rho t}|\nabla E(x)|^2dt=\frac{1}{2\rho}|\nabla E(x)|^2.
\end{multline*}
On a ainsi démontré l'inégalité suivante 
\begin{equation}
\label{eq-1}
E(x)-E(\beta)\leq \frac{1}{2\rho}|\nabla E(x)|^2=-\frac{1}{2\rho}\frac{d}{dt}E(S_t(x))\Big|_{t=0},
\end{equation}
qui est une inégalité entre l'entropie et la production d'entropie. Elle peut apparaître comme une simple inégalité de convexité mais elle a des conséquences intéressantes. 

On remarque que cette inégalité est optimale au sens où  pour l'exemple  classique, $E(x)=|x|^2/2$ on a une égalité. Par ailleurs, le flot de gradient $(S_t(x))_{t\geq 0}$ est une  interpolation remarquable entre $x$ et $\beta$ (lorsque $t=0$ et $t=\infty$), permettant de démontrer une inégalité de type~\eqref{eq-1}.

\medskip

$\bullet$ {\bf Convergence à l'équilibre du flot de gradient}

De cette inégalité~\eqref{eq-1} on exhibe  un taux explicite de la convergence à l'équilibre du flot de gradient. Nous savions déjà que le flot $(S_t(x))_{t\geq 0}$ convergeait vers $\beta$ mais on peut maintenant préciser la vitesse et l'espace naturel.   En effet,  il est facile de démontrer que l'inégalité entropie-production d'entropie~\eqref{eq-1} (pour un certain $\beta$) est  équivalente à la convergence exponentielle en entropie du flot de gradient. Plus précisément pour tout $x\in\R^n$ et $t\geq0$, 
$$
0\leq E(S_t(x))-E(\beta)\leq e^{-2\rho t}(E(x)-E(\beta)).
$$
Notons que l'équivalence n'est plus vérifiée si la convergence est de la forme  $Ce^{-2\rho t}$ avec $C>1$ au lieu de simplement  $e^{-2\rho t}$, dans ce cas on parle d'hypocoercivité et ces techniques doivent alors être modifiées. 

Cette simple  méthode a été utilisée un très grand nombre de fois depuis  1985 l'année de parution du papier de Bakry-\'Emery au {\it Séminaire   de probabilités},~\cite{bakry-emery1985}. Nous l'illustrerons avec deux exemples remarquables dans la section suivante. 

% Ainsi, dans ce modèle élémentaire, le cadre entropique est celui qui est le plus agréable pour l'étude du flot de gradient.

\section{Que se passe-t-il en dimension infinie ?}
\label{sec-3}

Une généralisation intéressante aux flots de gradient en dimension infinie a été proposée par Félix Otto en particulier dans~\cite{otto2001} et plus tôt dans l'algorithme JKO  proposé par Jordan, Kinderlehrer et Otto, cf.~\cite{jko1998}. C'est Villani qui lui a donné son nom, {\it le calcul d'Otto},  dans son ouvrage de référence~\cite[Chapitre 15]{villani2009}.

Ces travaux apportent une nouvelle utilisation de la distance de Wasserstein. Depuis les travaux de Kantorovich en 1942, cette distance sur l'espace des mesures de probabilités était largement utilisée en théorie des probabilités et en statistique pour estimer des convergences ou des déviations mais peu en EDP.   L'idée d'Otto est simple~: montrer qu'un flot de mesures comme par exemple l'équation de la chaleur est simplement le flot de gradient d'une fonctionnelle sur l'espace des mesures, considéré avec une métrique adaptée.  Considérer l'équation de la chaleur comme un flot de gradient n'est pas une idée nouvelle mais  Otto propose une nouvelle métrique, plus naturelle que nous détaillons ici.   

%\medskip

Nous nous permettons ici, de ne pas être rigoureux. De nombreux problèmes techniques ont été réglés en partie dans les travaux de Ambrosio, Gigli et Savaré, on pourra par exemple consulter l'ouvrage de référence~\cite{ambrosio-gigli2008}. Soit $M=\mathcal P_2(\R^n)$, l'ensemble des mesures de probabilités dans $\R^n$, absolument continues par rapport à la mesure de Lebesgue et ayant un moment d'ordre~2 fini. On confond de nouveau la mesure de probabilité avec sa densité par rapport à la mesure de Lebesgue.  

%\medskip

Otto a l'idée de considérer $M$ comme une variété riemannienne de dimension infinie où la distance de Wasserstein est simplement la distance riemannienne. 

\medskip

$\bullet$ {\bf Equation de continuité et vélocité} 

Soit un chemin dans l'ensemble des mesures de probabilités, 
$$
t\mapsto \mu_t\in M=\mathcal P_2(\R^n), 
$$
alors il existe une unique fonction $t\mapsto \Phi_t$ (à une constante près)   telle que 
\begin{equation}
\label{eq-11}
\partial_t\mu_t=-{\rm Div}(\mu_t \nabla \Phi_t),
\end{equation}
où ${\rm Div}$ est l'opérateur divergence dans $\R^n$. Bien entendu, pour que l'équation de continuité~\eqref{eq-11} soit vérifiée, le chemin $(\mu_t)_{t\geq0}$ doit vérifier quelques propriétés. Sans entrer dans les détails, il suffit que le chemin soit absolument continu dans l'espace de Wasserstein dans un sens expliqué dans~\cite[Chap.~8]{ambrosio-gigli2008}, et dans ce cas l'équation~\eqref{eq-11} est vérifiée au sens faible. 

\medskip

Nous allons identifier  la quantité $\partial_t\mu_t$ par la fonction $\nabla \Phi_t$ et on note  la vélocité du chemin par l'identification suivante
$$
\dot{\mu}_t=\nabla \Phi_t, 
$$
où on utilise comme en dimension finie la notation provenant de la mécanique. Cette représentation permet de voir l'évolution d'un flot de probabilités, où $\dot{\mu}_t$ est un vecteur qui montre la direction et la vitesse d'un élément de masse. Cette notion provient de la mécanique des fluides. 
Utilisée de cette façon, {\it la vélocité} (velocity) peut apparaître  un anglicisme plaisant, il  pourrait être remplacé par {\it  le vecteur vitesse}.

\medskip

$\bullet$ {\bf Espace tangent et métrique } 

Par l'équation de continuité, on a défini la vélocité du chemin $(\mu_t)_{t\geq 0}$ (ayant suffisamment de régularité). L'espace tangent est donc naturellement défini à partir de la vélocité,  soit $\mu\in M$, alors
$$
T_\mu M=\{\nabla\Phi, \Phi:\R^n\mapsto \R\}.
$$
De façon rigoureuse, l'espace tangent est l'adhérence de cet espace pour des fonctions régulières à support compact.  
La métrique proposée par Otto est la suivante, pour tout $\nabla\Phi,\nabla\Psi\in T_\mu M$, 
\begin{equation}
\label{eq-8}
\langle\nabla\Phi,\nabla\Psi\rangle_\mu=\int \nabla\Phi\cdot \nabla\Psi d\mu.
\end{equation}
\medskip

$\bullet$ {\bf Distance de Wasserstein}

La distance de Wasserstein, qui s'appelle aussi (plus justement) distance de Monge-Kantorovich,  est définie par~: soit $\mu,\nu\in M$,
$$
W_2(\mu,\nu)=\inf\sqrt{ \iint  |x-y|^2d\pi(x,y)},
$$ 
où l'infimum est pris sur l'ensemble des probabilités $\pi\in \mathcal  P(\dR^n\times \dR^n)$ ayant pour marginales $\mu$ et $\nu$. 

Dans~\cite{benamou-brenier2000}, J.-D. Benamou et Y. Brenier montrent que pour $\mu,\nu\in M$, mesures  régulières et ayant un moment d'ordre~2,  
$$
W_2(\mu,\nu)=\inf \sqrt{\int_0^1 \int |\nabla \Phi_t|^2d\mu_t dt},
$$ 
où l'infimum, pris sur l'ensemble des chemins $(\mu_t)_{t\in[0,1]}$ joignant $\mu$ à $\nu$ et $\nabla\Phi_t$, est la vélocité du chemin $(\mu_t)_{t\in[0,1]}$ considéré. 

En d'autres termes, la distance de Wasserstein est la distance riemannienne sur $M$ associée à la métrique définie précédemment en~\eqref{eq-8}. La métrique définie par le produit scalaire~\eqref{eq-8} est souvent  appelée métrique d'Otto.

\medskip

$\bullet$ {\bf Gradient et hessienne  d'une fonctionnelle}

De la  métrique, on peut  calculer le gradient d'une fonctionnelle, élément de l'espace tangent. Prenons par exemple  l'entropie de Boltzmann définie en~\eqref{eq-9}. Si on prend un chemin $(\mu_t)_{t\geq0}$ de vélocité $\dot \mu_t=\nabla\Phi_t$, alors 
$$
\frac{d}{dt}{\rm Ent}(\mu_t)=\frac{d}{dt}\int \mu_t\log\mu_t=-\int \log\mu_t{\rm Div}( \mu_t\nabla\Phi_t)=\int \nabla\log\mu_t\cdot\nabla\Phi_t\, \mu_t=\langle\nabla\log\mu_t,\dot\mu_t \rangle_{\mu_t}.
$$
Ainsi,  le gradient de l'entropie au point $\mu\in M$, que l'on note dans ce contexte  ${\rm grad}_\mu{\rm Ent}$, est donné par, 
$$
{\rm grad}_\mu{\rm Ent}=\nabla\log\mu.
$$

%\medskip

De la métrique, on peut aussi définir une dérivée covariante, basée sur une connexion riemannienne. Par ce biais, on peut définir les géodésiques (appelées géodésiques de McCann) dans l'espace de Wasserstein et pour finir la hessienne d'une fonctionnelle.  

En effet si $\mathcal F:M\mapsto \R$ est une fonctionnelle sur $M$, et si $(\mu_s)_{s\in[0,1]}$ une géodésique (à vitesse constante) de vélocité $(\dot\mu_s)_{s\in[0,1]}$, on peut définir la hessienne au point $\mu_s$, appliquée au vecteur tangent $\dot\mu_s$,  de la façon suivante~: 
$$
{\rm Hess}_{\mu_s}\mathcal F(\dot\mu_t,\dot\mu_t )=\frac{d^2}{ds^2}\mathcal F(\mu_s),
$$
si tout est assez régulier pour dériver deux fois. Pour obtenir une formule explicite nous avons besoin de l'équation des géodésiques.  Sans entrer dans les détails techniques, en particulier les problèmes délicats de régularité,  une géodésique $(\mu_s)_{s\in[0,1]}$ vérifie de façon formelle le système différentiel suivant 
\begin{equation}
\label{eq-12}
\left\{
\begin{array}{l}
\disp\partial_s\mu_s+{\rm Div}(\mu_s\nabla\Phi_s)=0,\\
\disp\partial_s\Phi_s+\frac{1}{2}|\nabla\Phi_s|^2=0.
\end{array}
\right.
\end{equation}
Par ce biais, au moins de façon formelle, on peut  calculer explicitement  $\frac{d^2}{ds^2}\mathcal F(\mu_s)$ et obtenir une expression de la hessienne.  

Reprenons l'exemple de l'entropie de Boltzmann. Soit $(\mu_s)_{s\in[0,1]}$ une géodésique de vélocité $(\nabla\Phi_s)_{s\in[0,1]}$,  on a vu précédemment que 
$$
\frac{d}{ds}{\rm Ent}(\mu_s)=\int \nabla\log\mu_s\cdot\nabla\Phi_s\, \mu_s=\int \nabla\mu_s\cdot\nabla\Phi_s.
$$
Ainsi, de façon formelle en supposant que tout est régulier, en particulier pour le système~\eqref{eq-12}, 
\begin{multline*}
\frac{d^2}{ds^2}{\rm Ent}(\mu_s)=\int \nabla\partial_s\mu_s\cdot\nabla\Phi_s+\int \nabla\mu_s\cdot\nabla\partial_s\Phi_s=\\
-\int \nabla{\rm Div}(\mu_s\nabla\Phi_s)\cdot\nabla\Phi_s-\frac12\int \nabla\mu_s\cdot\nabla|\nabla\Phi_s|^2=\\
\int \PAR{-\nabla\Phi_s\cdot\nabla\Delta\Phi_s+\frac12 \Delta|\nabla\Phi_s|^2}\mu_s,
\end{multline*}
où la dernière ligne s'obtient après 3 intégrations par parties.

Par ce calcul, au point $\mu\in M$ et appliquée au vecteur tangent $\nabla \Phi\in T_\mu M$, l'expression de la hessienne de l'entropie est donnée par   
$$
{\rm Hess}_\mu {\rm Ent} (\nabla\Phi,\nabla\Phi)=\int \PAR{\frac12\Delta|\nabla\Phi|^2-\nabla\Phi\cdot\nabla\Delta\Phi}\mu.
$$

On a donc maintenant tous les ingrédients pour appliquer en dimension infinie, ce qui a été fait précédemment en dimension finie.

\subsubsection*{Illustration avec deux exemples remarquables~: } 

Même si ce qui précède est formel, dû en particulier  aux problèmes de régularité, dans les deux exemples suivants, les résultats sont rigoureux car ils reposent sur des équations de types paraboliques qui ont des solutions régulières. 

\begin{enumerate}[\bf 1.]
\item {\bf \'Equation de la chaleur. }

Considérons  l'équation de la chaleur dans  $\R^n$, 
$$
\partial_t\mu_t=\Delta \mu_t={\rm Div}(\nabla \mu_t)={\rm Div}(\mu_t\nabla \log\mu_t).
$$
Ainsi, dans le calcul d'Otto, la vélocité du chemin de probabilité $\R^+\ni t\mapsto\mu_t$, est donnée par
$$
\dot{\mu}_t=-\nabla\log\mu_t.
$$
Si on reprend l'entropie de Boltzmann  ${\rm Ent}(\mu)=\int \mu\log\mu$,  on a vu précédemment que ${\rm grad}_\mu{\rm Ent}=\nabla\log\mu$,  on a 
$$
\dot{\mu}_t=-{\rm grad}_{\mu_t} {\rm Ent }.
$$
Ainsi l'équation de la chaleur est le flot de gradient de l'entropie de Boltzmann dans la métrique d'Otto.  Ce résultat est fondamental, démontré dans~\cite{jko1998}, article précurseur de ce domaine. Il permet aux auteurs de proposer un algorithme, basé sur la distance de Wasserstein, pour approcher la solution de l'équation de la chaleur.

Ainsi, tout ce qui a été démontré dans le cas fini dimensionnel s'applique. Il reste à 
comprendre quand la hessienne de la fonctionnelle $\rm Ent$ est uniformément minorée par une constante strictement positive. On a vu précédemment que  
$$
{\rm Hess}_\mu {\rm Ent} (\nabla\Phi,\nabla\Phi)=\int \PAR{\frac12\Delta|\nabla\Phi|^2-\nabla\Phi\cdot\nabla\Delta\Phi}\mu=\int ||{\rm Hess\, }\Phi||^2\mu,
$$
où $ ||{\rm Hess\, }\Phi||^2=\sum_{i,j=1}^n(\partial_{i,j}\Phi)^2$. Maintenant, si on souhaite obtenir une inégalité de type~\eqref{eq-2}, soit donc
$$
{\rm Hess}_\mu {\rm Ent} (\nabla\Phi,\nabla\Phi)=\int ||{\rm Hess\, }\Phi||^2\mu\geq \rho \int |\nabla\Phi|^2\mu, 
$$
pour toute fonction $\Phi$. Cette inégalité n'est malheureusement possible que si  $\rho=0$ (prendre par exemple $\Phi(x)=a\cdot x$, $a\in\R^n$). 

C'est assez naturel car,  dans $\R^n$, la solution de l'équation de la chaleur tend vers 0 (qui n'est pas une probabilité, de la masse a été perdue).  Pour éviter ça, on peut, soit se placer dans une variété riemannienne compacte, soit modifier l'équation de la chaleur et ajouter une force de rappel, lui permettant de converger vers une distribution non nulle. 
 
Ainsi, pour rester dans l'espace euclidien,  nous considérons l'équation de Fokker-Planck, 
$$
\partial_t\mu_t=\Delta \mu_t+{\rm Div}\PAR{x\mu_t}={\rm Div}\Big(\mu_t\nabla\SBRA{ \log\mu_t+\frac{|x|^2}{2}}\Big).
$$
Cette équation admet une solution explicite et régulière partant par exemple d'une donnée initiale dans l'espace de Schwartz. Elle admet une unique solution stationnaire de masse~1,  la mesure gaussienne standard, 
$$
d\gamma(x)=\exp\Big(\!-\frac{|x|^2}{2}\Big)\frac{1}{(2\pi)^{n/2}}d\lambda(x).
$$

La vélocité de la solution $(\mu_t)_{t\geq0}$ est par définition, 
$$
\dot{\mu}_t=-\nabla\Big(\log\mu_t+\frac{|x|^2}{2}\Big),
$$
et est le flot de gradient de la fonctionnelle 
$$
\mathcal F(\mu)=\int \mu\log\mu+\frac{|x|^2}{2}\mu,
$$
 par rapport à la métrique d'Otto. En effet, le gradient d'une fonctionnelle est une application linéaire et on obtient à la mesure de probabilité $\mu$, 
 $$
 {\rm grad}_\mu\mathcal F=\nabla\log \mu+ \nabla\frac{|x|^2}{2}.
 $$
 Il n'est pas difficile de montrer que la hessienne de $\mathcal F$ est uniformément minorée par constante 
 $\rho=1$, c'est-à-dire dans le langage d'Otto on peut montrer que  
 $$
 {\rm Hess}_\mu\mathcal F (\nabla\Phi,\nabla\Phi)=\int ||{\rm Hess\, }\Phi||^2\mu+\int {\rm Hess}\Big(\frac{|x|^2}{2}\Big)(\nabla \Phi,\nabla \Phi)\mu\geq \int |\nabla\Phi|^2\mu.
 $$
 De façon condensée, cette inégalité s'écrit simplement 
 $$
 {\rm Hess}_\mu\mathcal F \geq {\rm Id}_\mu,
 $$
 où $\rm{Id}_\mu$ représente la métrique au point $\mu\in M$. On retrouve la même inégalité qu'en dimension finie, l'inégalité~\eqref{eq-2} avec la constante $\rho=1$.

Le minimum de la fonctionnelle $\mathcal F$ est atteint pour la mesure gaussienne standard $\gamma$.  Ainsi, en imitant, ce qui se passe en dimension finie, la méthode de Bakry-\'Emery permet de démontrer l'inégalité suivante 
  $$
  \int \Big|\nabla\Big({\log \mu_t+\frac{|x|^2}{2}}\Big)\Big|^2\mu_t\leq e^{-2t} \int \Big|\nabla\Big(\log \mu_0+\frac{|x|^2}{2}\Big)\Big|^2\mu_0,
  $$
  inégalité équivalente à~\eqref{eq-10} dans le cas fini dimensionnel. 
L'inégalité entropie-production d'entropie s'écrit de la façon suivante, pour tout $\mu\in M$, 
 $$
 \mathcal F(\mu)-\mathcal F(\gamma)\leq \frac{1}{2}|{\rm grad}_\mu\mathcal F|^2_\mu=\frac{1}{2}\int \Big|\nabla\Big(\log \mu+\frac{|x|^2}{2}\Big)\Big|^2\mu. 
 $$
 Après un changement de fonction, $\mu=\exp(-|x|^2/2)f$, cette inégalité prend la forme suivante, pour toute fonction $f$ positive,  
 $$
 \int f\log \frac{f}{\int fd\gamma}d\gamma\leq \frac{1}{2}\int \frac{|\nabla f|^2}{f}d\gamma.
 $$
Celle-ci s'appelle l'inégalité de Sobolev logarithmique pour la mesure gaussienne $\gamma$ et a été démontrée par Gross en 1975 dans~\cite{gross}. Cette inégalité a de nombreuses propriétés remarquables et est utilisée en EDP, probabilités, géométrie etc. 

\item {\bf L'équation des milieux poreux avec une force de rappel. }

Considérons maintenant l'équation de Fokker-Planck non linéaire dans $\R^n$, $n>2$ , ou équation des milieux poreux avec une dérive (en fait plutôt appelée équation des diffusions rapides puisque l'exposant, ici $1-1/n$, est  dans l'intervalle $]0,1[$),  
\begin{equation}
\label{eq-7}
\partial_t\mu=\frac{1}{n}\Delta\PAR{\mu^{1-1/n}}+\frac{n-1}{n}{\rm Div}\PAR{\mu x}=\frac{n-1}{n}{\rm Div}\SBRA{\mu\nabla\PAR{-\mu^{-1/n}+\frac{|x|^2}{2}}}.
\end{equation}
Partant d'une donnée initiale qui est une probabilité, la vélocité du flot   est donc donnée par 
$$
\dot{\mu}_t=-\frac{n-1}{n}\nabla \PAR{-\mu_t^{1/n}+\frac{|x|^2}{2}}.
$$
Il n'est pas trop difficile de voir que 
$$
\dot{\mu}_t=-{\rm grad}_{\mu_t} {\rm \mathcal F },
$$
où
$$
\mathcal F (\mu)=\int\PAR{ -\mu^{-1/n}+\frac{n-1}{n}\frac{|x|^2}{2}}\mu.
$$
Un calcul un peu plus compliqué, si on ne s'y prend mal, permet de calculer la hessienne de $\mathcal F$. Au point $\mu$ et pour $\nabla\Phi$, élément de l'espace tangent,  on obtient 
$$
{\rm Hess}_\mu\mathcal F (\nabla\Phi,\nabla\Phi)=\frac1n\int (||{\rm Hess\, }\Phi||^2-\frac1n(\Delta\Phi)^2)\mu+\frac{n-1}{n}\int |\nabla\Phi|^2\mu,
$$
calcul effectué par exemple dans~\cite{otto-westdickenberg2006}. 
Puisque $||{\rm Hess\, }\Phi||^2=\sum_{i,j}(\partial_{ij}\Phi)^2$, par une inégalité de Cauchy-Schwarz,  
$$
||{\rm Hess\, }\Phi||^2-\frac1n(\Delta\Phi)^2\geq0.
$$
Ainsi on obtient,  
$$
{\rm Hess}_\mu\mathcal F (\nabla\Phi,\nabla\Phi)\geq \frac{n-1}{n}\int |\nabla\Phi|^2\mu,
$$
ce qui s'écrit aussi de la façon suivante,
$$
{\rm Hess}_\mu\mathcal F \geq \frac{n-1}{n}{\rm Id_\mu}.
$$
On retrouve encore le cas de la dimension finie~\eqref{eq-2} avec la constante $\rho=\frac{n-1}{n}$. Un calcul indépendant montre que la solution 
$(\mu_t)_{t\geq0}$ converge, lorsque $t$ tend vers l'infini, vers l'état stationnaire  $\mu_\infty=(C+|x|^2/2)^{-n}$ où $C$ est une constante de normalisation pour que $\mu_\infty$ soit une mesure de probabilité. 
On peut alors  montrer une inégalité d'entropie-production d'entropie, pour toute mesure $\mu_0\in M$, 
$$
\mathcal F(\mu_0)-\mathcal F(\mu_\infty)\leq\frac{n-1}{2n}\int \left|\nabla\PAR{\mu_0^{-1/n}-\frac{|x|^2}{2}}\right|^2\mu_0=- \frac{n}{2(n-1)}\frac{d}{dt}\mathcal F(\mu_t)\Big|_{t=0}, 
$$
Cette inégalité de convexité n'est rien d'autre qu'une forme un peu «~barbare~», et équivalente, de l'inégalité de Sobolev optimale dans $\R^n$, pour toute fonction $f:\R^n\mapsto \R$, régulière à support compact, 
$$
||f||_{L^{\frac{2n}{n-2}}(\lambda)}\leq C_{op}||\nabla f||_{L^2(\lambda)}.
$$
Cette dernière inégalité se retrouve simplement après un changement de fonction et une intégration par parties. Notons que cette inégalité est attribuée à Sobolev en 1938 et que la constante optimale a été calculée par Aubin et Talenti en 1976. 

Il est remarquable de voir qu'une telle inégalité, si importante en analyse, en EDP et en géométrie, n'est qu'une inégalité de convexité le long d'un flot de gradient. On s'aperçoit aisément dans la preuve de cette inégalité que d'autres cas assez similaires peuvent être traités, par exemple on retrouve naturellement par cette méthode la famille des inégalités de Gagliardo-Nirenberg optimales démontrées par M. del Pino et J. Dolbeault en 2002 ou bien l'inégalité de Sobolev dans une variété riemannienne avec des conditions de courbures, démontrée par J. Demange en 2008. 
\end{enumerate}

L'interaction entre l'entropie et le calcul d'Otto est toujours un domaine très actif en mathématique. Pour l'illustrer, écrivons en détail ce récent théorème de S. Zugmeyer qui démontre une inégalité de type Sobolev  sur un ensemble convexe borné de $\R^n$,  mettant  à profit le calcul d'Otto pour généraliser les résultats précédents. 

\medskip 

\noindent 
{\bf Théorème (\cite{zug})} {\it  Soit $H:\R^+\mapsto \R$ une  fonction de classe $\mathcal C^2$, strictement convexe, vérifiant $H(0)=0$ et posons $\Psi=H'$.  Choisissons un ensemble convexe fermé $\Omega\subset\R^n$ ($n\geq1$) et une fonction de classe $\mathcal C^2$, $v:\Omega\mapsto (0,+\infty)$. Supposons que 
\begin{itemize}
\item Pour tout $x\geq0$, $xU'(x)+\frac{1-n}{n}U(x)\geq0$, où $U(x)=x\Psi(x)-H(x)$; 
\item Il existe une constante $C>0$, telle que $-{\rm Hess}(\Psi(v))\geq C{\rm Id}$.
\end{itemize}
Alors pour tout fonction strictement positive $u\in\mathcal C^\infty(\Omega)$ vérifiant $\int_\Omega u=\int_\Omega v$,  
$$
\int_{\Omega}\PAR{H(u)-H(v)-(u-v)\Psi(v)}d\lambda\leq\frac{1}{2C}\int_\Omega ||\nabla\PAR{\Psi(v)-\Psi(u)}||^2ud\lambda.
$$
}

%{\bf Conclusion :} Clausius a pu inventer un mot qui est entré dans le langage courant, mot désignant parfois le désordre. Par ailleurs, autour de ce mot, on trouve une notion ou bien un concept assez flou permettant son utilisation dans des domaines mathématiques contemporains. L'exemple de l'inégalité optimale de Sobolev est, il me semble, un exemple frappant. 

\bigskip

Terminons ces quelques pages par une citation de Ievgueni Zamiatine, tiré de {\it Nous autres} (1920), \og Alors voici~: il y a deux forces en ce monde - l'entropie et l'énergie. L'entropie vise la paix et la béatitude, l'équilibre heureux - l'autre recherche la rupture des équilibres, la torture du mouvement infini.\fg
C'est très probablement la première utilisation de l'entropie en littérature. On peut se demander si l'auteur de ce remarquable livre d'anticipation reprend les mêmes notions de l'entropie que celles proposées par Clausius. Si Clausius n'a pas pu donner son avis sur cette question, le régime soviétique lui n'a pas apprécié cette œuvre satirique  dénonçant le totalitarisme. 

\medskip

Je tiens à remercier Louis Dupaigne et Laurent Miclo pour leurs remarques pertinentes lors de la rédaction de ces quelques pages.

\medskip

\small

\noindent
I. G. Institut Camille Jordan, Umr Cnrs 5208, Universit\'e Claude Bernard Lyon 1, 43 boulevard du 11 novembre 1918, F-69622 Villeurbanne cedex.
    \texttt{gentil@math.univ-lyon1.fr}

%\bibliographystyle{alpha}
%{\footnotesize{\bibliography{biblio}}}

\begin{thebibliography}{{Gro}75}

\bibitem[AGS08]{ambrosio-gigli2008}
L.~{Ambrosio}, N.~{Gigli}, and G.~{Savar\'e}.
\newblock {\em {Gradient flows in metric spaces and in the space of probability
  measures.}}
\newblock Basel: Birkh\"auser, 2nd edition, 2008.

\bibitem[BB00]{benamou-brenier2000}
J.-D. {Benamou} and Y.~{Brenier}.
\newblock {A computational fluid mechanics solution to the Monge-Kantorovich
  mass transfer problem.}
\newblock {\em {Numer. Math.}}, 84(3):375--393, 2000.

\bibitem[BE85]{bakry-emery1985}
D.~{Bakry} and M.~{\'Emery}.
\newblock {Diffusions hypercontractives.}
\newblock {\em {S\'emin. de probabilit\'es XIX, Univ. Strasbourg 1983/84,
  Proc., Lect. Notes Math. 1123, 177-206}}, 1985.

\bibitem[Cla65]{clausius}
R.~Clausius.
\newblock Sur diverses formes facilement applicables qu’on peut donner aux
  {\'e}quations fondamentales de la th{\'e}orie m{\'e}canique de la chaleur.
\newblock {\em {J. Math. Pures Appl. (2)}}, 10:361--400, 1865.

\bibitem[{Gro}75]{gross}
L.~{Gross}.
\newblock {Logarithmic Sobolev inequalities.}
\newblock {\em {Am. J. Math.}}, 97:1061--1083, 1975.

\bibitem[JKO98]{jko1998}
R.~{Jordan}, D.~{Kinderlehrer}, and F.~{Otto}.
\newblock {The variational formulation of the Fokker-Planck equation.}
\newblock {\em {SIAM J. Math. Anal.}}, 29(1):1--17, 1998.

\bibitem[{Ott}01]{otto2001}
F.~{Otto}.
\newblock {The geometry of dissipative evolution equations: The porous medium
  equation.}
\newblock {\em {Commun. Partial Differ. Equations}}, 26(1-2):101--174, 2001.

\bibitem[OW06]{otto-westdickenberg2006}
F.~{Otto} and M.~{Westdickenberg}.
\newblock {Eulerian calculus for the contraction in the Wasserstein distance.}
\newblock {\em {SIAM J. Math. Anal.}}, 37(4):1227--1255, 2006.

\bibitem[{Per}02]{perelman}
G.~{Perelman}.
\newblock {The entropy formula for the Ricci flow and its geometric
  applications}.
\newblock {\em Preprint}, 2002.

\bibitem[{Vil}09]{villani2009}
C.~{Villani}.
\newblock {\em {Optimal transport. Old and new.}}
\newblock Berlin: Springer, 2009.

\bibitem[Zug19]{zug}
S.~Zugmeyer.
\newblock {Entropy flows and functional inequalities in convex sets}.
\newblock Preprint hal-02407995, 2019.

\end{thebibliography}

\end{document}